\documentstyle[12pt]{article}
\title{Decomposition of $q$-deformed Fock spaces}
\author{Masaki Kashiwara, Tetsuji Miwa, and Eugene Stern}
\date{Draft, August 1995}

\begin{document}
\newtheorem{thm}{Theorem}[section]
\newtheorem{prop}{Proposition}[section]
\newtheorem{lemma}{Lemma}[section]
\newtheorem{rem}{Remark}[section]

\font\germ=eufm10
\def\goth#1{\hbox{\germ#1}}
\def\spc{\phantom{l}}
\def\be{\begin{equation}}
\def\exp{{\rm exp}}
\def\en{\end{equation}}
\def\bea{\begin{eqnarray}}
\def\ena{\end{eqnarray}}
\def\bean{\begin{eqnarray*}}
\def\enan{\end{eqnarray*}}
\def\lb#1{\label{eqn:#1}}
\def\rf#1{(\ref{eqn:#1})}
\def\qed{\hfill\fbox{}\break\smallskip}
\def\pf{\noindent{\it Proof.\quad}}
\def\displ#1{{\displaystyle #1}}
\def\XXX{{X\hskip-2pt X\hskip-2pt X}}
\def\XXZ{{X\hskip-2pt X\hskip-2pt Z}}
\def\XYZ{{X\hskip-2pt Y\hskip-2pt Z}}
\def\mod{{\rm mod}\,}
\def\slh{\widehat{\goth{sl}}}
\def\glh{\widehat{\goth{gl}}}
\def\slth{\slh_2}
\def\slnh{\slh_n}
\def\glth{\glh_2}
\def\glnh{\glh_n}
\def\op{\scriptstyle {\rm op}}
\def\V{{\cal V}}
\def\uq{U_q(\slth)}
\def\uqn{U_q(\slnh)}
\def\upqn{U'_q(\slnh)}
\def\uqc{U_q(\slt)}
\def\uqglt{U_q(\glth)}
\def\End{{\rm End}}
\def\Hom{{\rm Hom}}
\def\Ker{{\rm Ker}}
\def\Im{{\rm Im}}
\def\id{{\rm id}\,}
\def\ch{{\rm ch}\,}
\def\tr{{\rm tr}\,}
\def\ad{{\rm ad}\,}
\def\qdet{\hbox{$q$-det}}
\def\ket#1{|#1\rangle}
\def\tPhi{{\tilde\Phi}}
\def\vaf{V_{\hbox{\scriptsize aff}}}
\def\F{{\cal F}}
\def\H{{\cal H}}
\def\Z{{\bf Z}}
\def\Q{{\bf Q}}
\def\C{{\bf C}}
\def\e{\varepsilon}
\def\vep{\varepsilon}
\def\z{\zeta}
\def\cL{{\cal L}}
\def\tL{\widetilde{\cal L}}
\def\bL{\overline{L}}
\def\cR{{\cal R}}
\def\bR{\overline{R}}
\def\cM{{\cal M}}
\def\cN{{\cal N}}
\def\tN{\widetilde{\cN}}
\def\hV{\widehat{\cal V}}
\def\Vh{\widehat{\cal V}}

\def\slhat{\slnh}
\def\qslhat{\upqn}
\def\qslhatext{\uqn}
\def\hNhat{\widehat{H}_N(q^2)}
\def\hinfhat{\widehat{H}_{\infty}(q^2)}
\def\vac#1{|#1\rangle}



\maketitle
\begin{abstract}
A decomposition of the level-one $q$-deformed Fock
representations of $\uqn$ is given. It is found that
the action of $\qslhat$ on these Fock spaces is centralized by a
Heisenberg algebra,
which arises from the center of the affine
Hecke algebra $\widehat{H}_N$ in the limit $N\rightarrow\infty$.
The $q$-deformed Fock space
is shown to be isomorphic as a
$\upqn$-Heisenberg-bimodule to the tensor product
of a level-one irreducible highest weight representation of
$\qslhat$ and the Fock representation of the Heisenberg algebra.
The isomorphism is used to decompose the $q$-wedging operators,
which are intertwiners between the $q$-deformed Fock spaces,
into constituents coming from $\qslhat$ and from the Heisenberg algebra.
\end{abstract}

\section{$q$-wedges and $q$-deformed Fock space}
This first introductory section describes the realization of $q$-deformed
Fock space in terms of $q$-wedges.  This Fock space
was first constructed by a different method by Hayashi in \cite{Hayashi},
and described in terms of colored Young diagrams in \cite{Misra-Miwa}.
A less formal version of the $q$-wedge construction
was given in \cite{Eug}, out of which the exposition in this
section first evolved.

\subsection{Preliminaries on $\qslhatext$}
The algebras $\qslhatext$ and $\qslhat$ (\cite{Drinfeld,Jimbo})
will act on our Fock
space.  In this section, we will mainly work with $\qslhat$,
which is an algebra generated by elements $E_i$, $F_i$, and
$K_i^{\pm 1}$, $i=0,1,\ldots , n-1$, with the following relations if
$n>2$:
\begin{eqnarray}
&&K_i K_j = K_j K_i, \label{begingen} \\
&&K_i E_j= q^{a_{ij}} E_j K_i, \label{cartan-1} \\
&&K_i F_j = q^{-a_{ij}} F_j K_i, \label{cartan-2} \\
&&E_i F_j - F_j E_i = \delta_{i,j} \frac{K_i - K_i^{-1}}
  {q - q^{-1}}, \label{commutator} \\
&&E_i E_j = E_j E_i \hspace{100pt} \mbox{\rm if}
          \hspace{5pt} i\not=j\pm1, \\
&&F_i F_j = F_j F_i \hspace{105pt} \mbox{\rm if}
          \hspace{5pt} i\not=j\pm1, \\
&& E_i^2 E_{i \pm 1} - (q + q^{-1})E_i E_{i \pm 1} E_i
        + E_{i \pm 1} E_i^2=0, \\
&& F_i^2 F_{i \pm 1} - (q + q^{-1})F_i F_{i \pm 1} F_i
        + F_{i \pm 1} F_i^2=0. \label{endgen}
\end{eqnarray}
The indices in all these relations are to
be read modulo $n$.  In (\ref{cartan-1})--(\ref{cartan-2}),
$a_{ij}$ is $2$ if $i=j$, $-1$ if $i=j \pm 1$, and $0$ otherwise.

In the case $n=2$, we take $a_{ij}=-2$ if $i \neq j$.
Moreover, the last two relations (sometimes called
the $q$-Serre relations) are replaced by the following ones:
\begin{eqnarray}
E_i^3 E_{i \pm 1} - [3] E_i^2 E_{i \pm 1} E_i + [3] E_i
E_{i \pm 1} E_i^2 - E_{i \pm 1} E_i^3&=&0 ,\\
F_i^3 F_{i \pm 1} - [3] F_i^2 F_{i \pm 1} F_i + [3] F_i
F_{i \pm 1} F_i^2 - F_{i \pm 1} F_i^3&=&0 .
\end{eqnarray}
Here we have used the standard notation
$$[n] = \frac{q^n - q^{-n}}{q-q^{-1}}.$$
Throughout this paper, $q$ should be taken to be either a
formal parameter or a generic complex number (specifically,
not a root of unity).

$\qslhat$ is a Hopf algebra, with coproduct given by
\begin{eqnarray}
\Delta(K_i) &=& K_i \otimes K_i, \label{coprod-1} \\
\Delta(E_i) &=& E_i \otimes K_i + 1 \otimes E_i, \label{coprod-2} \\
\Delta(F_i) &=& F_i \otimes 1 + K_i^{-1} \otimes F_i . \label{coprod-3}
\end{eqnarray}

In the next section we will sometimes work with the extended quantum
affine algebra $\qslhatext$, which is a $1$-dimensional extension of
$\qslhat$ by the ``degree operator'' $d$, satisfying the relations
\begin{equation}
[d, K_i] = 0, \hspace{.8in} [d,E_i]=\delta_{i,0} E_i, \hspace{.8in}
[d,F_i] = -\delta_{i,0} F_i.
\label{DREL}
\end{equation}
The Hopf algebra structure extends naturally to $\qslhatext$ by
defining
\begin{equation}
\Delta(d) = d \otimes 1 + 1 \otimes d.
\end{equation}


The philosophy guiding this subject is that almost
any representation-theoretic construction involving
the affine Lie algebra $\slhat$ should have
an appropriate $q$-analog in the representation theory of $\qslhat$.
The aim of this section is to deform the infinite wedge (or ``Fermionic
Fock space'') level 1 representations of $\slhat$ to representations
of $\qslhat$.  (A very friendly introduction to $\slhat$
and its infinite
wedge highest weight representations can be found in \cite{Kac-Raina}.)

We will construct infinite $q$-wedges in terms of an infinite
tensor product of {\em evaluation modules}.  Let $V=\mbox{\bf C}^n$,
with basis $v_1, \ldots , v_n$,
and let $V(z) = V \otimes \mbox{\bf C}[z,z^{-1}]$, with basis
$\{z^av_j\}$.  Here $a \in \mbox{\bf Z}$ and $j=1,2,\ldots , n$,
so $V(z)$ is regarded as an infinite dimensional space.
$\uqn$ acts on $V(z)$ in the following way:
\begin{eqnarray}
K_i  (z^a v_j) &=& q^{\delta_{i,j} - \delta_{i+1,j}}
z^a v_j, \\
E_i  (z^a v_j) &=& \delta_{i,j-1}  z^{a+{\delta_{i,0}}}
 v_{j-1}, \\
F_i  (z^a v_j) &=& \delta_{i,j} z^{a-\delta_{i,0} }
 v_{j+1},\\
d(z^av_j) &=&a z^av_j .
\end{eqnarray}
The indices in all of these relations should be read modulo $n$.

The action of the $E_i$'s gives rise to a
a natural ordering on the basis $\{z^a v_j\}$.  This ordering is given by
\begin{equation}
\cdots > z^{a-1} v_2 > z^{a-1} v_1 > z^a v_n > z^a v_{n-1} > \cdots
> z^a v_1 > z^{a+1} v_n > \cdots \ , \label{ordering}
\end{equation}
and at times it will be convenient to relabel the basis in accordance
with it.  Namely, let $u_{j-an}= z^a v_j$; then $u_l > u_m$ just in
case $l>m$.  {\em In the course of working with $V(z)$ and its tensor
products, we will
sometimes use the $z^a v_j$ notation and sometimes the $u_m$ notation
(whichever happens to be more convenient) without further comment.}
The action on the basis $\{u_m\}$ is
\begin{eqnarray}
K_i  (u_m)&=&
q^{\delta(m\equiv i\ \mod\,n)-\delta(m\equiv i+1\ \mod\,n)}\,u_m,\\
E_i  (u_m) &=& \delta(m-1\equiv i\,\mod\ n)\,u_{m-1},\\
F_i  (u_m) &=& \delta(m\equiv i\,\mod\ n)\,u_{m+1}.
\end{eqnarray}
Here, for a statement $P$, $\delta(P)$ is equal to $1$
if $P$ is true and $0$ otherwise.
\par
Iterating the coproduct (\ref{coprod-1})-(\ref{DREL})
$N-1$ times defines a natural
action of $\uqn$ on the tensor product $V(z)^{\otimes N}$.
According to a ``quantum affine analog'' of the usual Weyl
duality between $GL_n$ and the symmetric group,
the centralizer of the action of $\qslhat \subset \uqn$
on $V(z)^{\otimes N}$ is the affine Hecke algebra
$\hNhat$.  (See \cite{GRV} for more
details.)  We turn next to this algebra.

\subsection{Preliminaries on the affine Hecke algebra}
$\hNhat$ is an associative algebra generated by elements
$T_i$, $i=1,\ldots ,N-1$, and $y_j^{\pm 1}$, $j=1,\ldots , N$.
These elements satisfy the following relations:
\begin{eqnarray}
T_i^2 &=& (q^2-1) \ T_i + q^2, \label{heckegen-1} \\
T_i T_{i+1} T_i &=& T_{i+1} T_i T_{i+1}, \label{heckegen-2} \\
T_i T_j &=& T_j T_i  \hspace{.8in} \mbox{\rm if} \hspace{2mm} |i-j|>1,
\label{heckegen-3} \\
y_i y_j &=& y_j y_i, \label{heckegen-4} \\
y_j T_i &=& T_i y_j \hspace{.8in} \mbox{\rm if $i \neq j,j+1$},
\label{heckegen-5} \\
y_j T_j &=& T_j y_{j+1} - (q^2-1) y_{j+1}, \label{heckegen-6} \\
y_{j+1} T_j &=& T_j y_j + (q^2 - 1) y_{j+1} . \label{heckegen-7}
\end{eqnarray}
Note that the $T_i$ are invertible because of (\ref{heckegen-1}),
with $T_i^{-1} = q^{-2} T_i + (q^{-2} - 1)$.  In light of this,
relations (\ref{heckegen-6}) and (\ref{heckegen-7})
are both equivalent to the relation
$$T_i y_i T_i = q^2 y_{i+1}.$$
Also, the relation (\ref{heckegen-1}) can be written as
\be
(T_i+1)(T_i-q^2)=0. \label{eigenvalues}
\en
Decomposing spaces on which the $T_i$
act into a $-1$-eigenspace and a $q^2$-eigenspace will be an important tool
in what follows.

The subalgebra $H_N(q^2) \subset \hNhat$ generated by just the $T_i$
is the usual (finite) Hecke algebra of type $A$, which is a $q$-deformation
of the symmetric group $S_N$.  The elements $T_i$ are the $q$-analogs of
the adjacent transpositions $\sigma_i = (i,\, i+1)$ in $S_N$.  In the same
way, $\hNhat$ should be thought of as a $q$-deformation of the {\em affine}
symmetric group, which is the Weyl group of the affine Lie algebra
$\slh_N$.

$\hNhat$ acts on $V(z)^{\otimes N}$
on the right in the following way.  First,
$y_j$ acts as multiplication by $z_j^{-1}$.  (Having it act by $z_j^{-1}$
rather than $z_j$ is necessary for compatibility
(\ref{heckegen-6})--(\ref{heckegen-7}) with the action
of $T_i$ defined below.)  The action of $T_i$
is given as follows.  Write elements $z^{a_1} v_{m_1}\otimes \cdots
\otimes z^{a_N} v_{m_N}$ as $(v_{m_1} \otimes \cdots \otimes v_{m_N})
\cdot z_1^{a_1} \cdots z_N^{a_N}$.  (This is the notation used in
\cite{GRV}.)  The symmetric group $S_N$ acts (by
permuting factors and variables, respectively) on both the tensor part
and the polynomial part of such an expression.  Write $(v_{m_1} \otimes
\cdots \otimes v_{m_N})^{\sigma_i}$ and $(z_1^{a_1} \cdots z_N^{a_N})^{
\sigma_i}$ for what results when $\sigma_i = (i,\,i+1)$ acts on
$v_{m_1} \otimes \cdots \otimes v_{m_N}$ and on $z_1^{a_1} \cdots
z_N^{a_N}$.  In terms of this
notation, the action of $T_i$ on $v_{m_1} \otimes \cdots \otimes v_{m_N}
\cdot \mbox{\bf z}$ (here $\mbox{\bf z}$ is shorthand for $z_1^{a_1}
\cdots z_N^{a_N}$) is given by the following set of formulas,
which are variants of the ones given in \cite{GRV}:
\begin{equation}
((v_{m_1} \otimes \cdots \otimes v_{m_N}) \cdot \mbox{\bf z})
\cdot T_i =
\left\{
\begin{array}{ll}
- q(v_{m_1} \otimes \cdots \otimes v_{m_N})^{\sigma_i} \cdot
\mbox{\bf z}^{\sigma_i} \\
\hspace{3mm} - (q^2-1)(v_{m_1} \otimes \cdots \otimes v_{m_N}) \cdot
\frac{z_{i+1} \mbox{\footnotesize \bf z}^{\sigma_i} -
z_i \mbox{\footnotesize \bf z} }{z_i-z_{i+1}} & \hspace{3mm}
\mbox{if $m_i<m_{i+1}$,} \\
\noalign{\smallskip}
- (v_{m_1} \otimes \cdots \otimes v_{m_N}) \cdot \mbox{\bf z}^{\sigma_i} \\
\hspace{3mm} - (q^2-1) (v_{m_1} \otimes \cdots \otimes v_{m_N}) \cdot
\frac{z_i(\mbox{\footnotesize \bf z}^{\sigma_i} -
\mbox{\footnotesize \bf z})}{z_i-z_{i+1}} & \hspace{3mm}
\mbox{if $m_i=m_{i+1}$,} \\
\noalign{\smallskip}
- q(v_{m_1} \otimes \cdots \otimes v_{m_N})^{\sigma_i} \cdot
\mbox{\bf z}^{\sigma_i} \\
\hspace{3mm} - (q^2-1)(v_{m_1} \otimes \cdots \otimes v_{m_N}) \cdot
\frac{z_i(\mbox{\footnotesize \bf z}^{\sigma_i} -
\mbox{\footnotesize \bf z})}{z_i-z_{i+1}} & \hspace{3mm}
\mbox{if $m_i>m_{i+1}$.}
\end{array} \right. \label{hecke-action}
\end{equation}
The important fact for us is that this right action of $\hNhat$
on $V(z)^{\otimes N}$ commutes with the left action
of $\qslhat$ on $V(z)^{\otimes N}$ given in the
previous subsection in terms of the coproduct.

\subsection{$q$-antisymmetrization and $q$-wedges}
As a vector space,
$H_N(q^2) \subset \hNhat$ has a natural basis $\{T_{\sigma}
\}_{\sigma\in S_N}$ which can be defined as follows.
Given $\sigma \in S_N$, take a minimal-length expression $\sigma =
\sigma_{i_1} \sigma_{i_2} \cdots \sigma_{i_l}$ of $\sigma$ in terms of
adjacent transpositions.  Then define $T_{\sigma} = T_{i_1} T_{i_2} \cdots
T_{i_l}$.  It is a basic result about the Hecke algebra that $T_{\sigma}$
depends only on $\sigma$ and not on the particular expression that was used.

We define
the $q$-antisymmetrizing operator $A^{(N)}$ acting on $V(z)^{\otimes N}$
to be the sum
\bea
A^{(N)} = \sum_{\sigma \in S_N}T_{\sigma}
\label{q-anti-symmetrizer}
\ena
There are no $(-1)^{\sigma}$ factors appearing in this sum
because the sign of the permutation is already incorporated
into the definition (\ref{hecke-action}) of the action of $\hNhat$;
for example, $(v_{i+1} \otimes v_i) \cdot T_1 = -q v_i \otimes v_{i+1}$.
(In other words, the operator $T_i$ acting in $V(z)^{\otimes N}$
is really a deformation of $-\sigma_i$, rather
than of $\sigma_i\in S_N$.)


The {\em $q$-antisymmetrization} of a pure tensor
$z^{a_1} v_{m_1} \otimes \cdots \otimes z^{a_N} v_{m_N}$
is defined to be
\begin{equation}
(z^{a_1} v_{m_1} \otimes \cdots \otimes z^{a_N} v_{m_N}) \cdot A^{(N)}.
\label{q-anti-sym}
\end{equation}

We begin our study of $A^{(N)}$ with the following proposition, which
essentially asserts that $A^{(N)}$ is (up to scalar) an idempotent,
and $V(z)^{\otimes N}$ decomposes as a direct sum of its two eigenspaces.

\begin{prop}\label{TSPLT}
\bea
V(z)^{\otimes N}=\Im A^{(N)}\oplus\Ker A^{(N)}.
\lb{AST1}
\ena
\end{prop}

\noindent
{\em Proof.} \hspace{2mm}
Note that for each $i=1,2,\cdots,N-1$ we have a factorization
\begin{equation}
A^{(N)}
=\left( \sum_{\sigma '} T_{\sigma '} \right)(T_i+1),
\label{facto}
\end{equation}
where $\sigma'$ ranges over $S_N / \{ \mbox{\rm id}, \sigma_i \}$.
{}From this and (\ref{eigenvalues}), it follows that
\[
A^{(N)}(T_i-q^2)=0.
\]
This means that the action of $T_i$ on the right on $\Im A^{(N)}$ is simply
multiplication by $q^2$.  Hence, right multiplication by $A^{(N)}$ on
$\Im A^{(N)}$ is equal to multiplication by the scalar
\begin{equation}
\sum_l n(l) q^{2l}, \label{sum-length}
\end{equation}
where $n(l)$ is the number of elements of $S_N$ having length $l$.
This sum is equal to the product
$$
\prod_{m=1}^N \frac{1-q^{2m}}{1-q^2},
$$
which is a non-zero scalar since $q$ is not a root of unity.
Therefore we have
\begin{equation}
A^{(N)}\left(A^{(N)}-\prod_{m=1}^N{1-q^{2m}\over1-q^2}\right)=0,
\end{equation}
and this implies the assertion.
\qed

In the classical ($q=1$) case, there are two equivalent ways of defining
the wedge product.  One is as the subspace of the tensor product consisting
of completely antisymmetric tensors (i.e., the image of the antisymmetrizer),
and the other is as a quotient of
the tensor product by relations of the form $v \wedge w = - w \wedge v$,
which generate the kernel of the antisymmetrizer.

In the quantum case ($q\not=1$),
both approaches are again available because of
Proposition \ref{TSPLT}.
The $q$-wedge space can be defined either as $\Im A^{(N)}$, the subspace
of $V(z)^{\otimes N}$ consisting of
$q$-antisymmetrized tensors, or as a quotient of this tensor product by
certain relations which generate the kernel of $A^{(N)}$.
Let us now describe
these relations.  The first step is

\begin{prop}
The kernel of $A^{(N)}$ is the sum of the kernels of the operators $T_i+1$,
$i=1,2,\ldots , N-1$. \label{A->adjrel}
\end{prop}

\noindent
{\em Proof.} \hspace{2mm}
Equation (\ref{facto}) implies that
$\sum_i \Ker (T_i+1) \subseteq \Ker A^{(N)}$.
Let us show that
$\Ker A^{(N)}\subseteq\sum_i\Ker(T_i+1)$.  Proposition
\ref{TSPLT} applied to $V(z)^{\otimes 2}$ asserts that
\begin{equation}
V(z)\otimes V(z)=\Ker(T-q^2)\oplus\Ker(T+1).
\label{SPLIT}
\end{equation}
In general, by arguing as we did for that Proposition, we can conclude
\begin{equation}
\omega A^{(N)} \equiv \left(\prod_{m=1}^N {1-q^{2m}\over1-q^2}\right)
\omega\bmod \sum_i\Ker(T_i+1) \label{prod-mod}
\end{equation}
for any $\omega\in V(z)^{\otimes N}$. This equation implies that if
$\omega\in\Ker A^{(N)}$, then $\omega\equiv0\bmod \sum_i\Ker(T_i+1)$.
\qed

Proposition \ref{A->adjrel}
shows that to find relations generating the kernel of
$A^{(N)}$, it suffices to find relations generating the kernel of each
$T_i+1$.  For ease of notation, let us restrict to considering $V(z)
\otimes V(z)$, on which $T=T_1$ acts.  Since $T$ commutes with the
action of $\qslhat$, $\Ker (T+1)$ is preserved by the action of $\qslhat$.
One way of deriving relations is to start with a simple element
of the kernel, say $v_1 \otimes v_1$, and then act on it by $\qslhat$.
This gives us the following elements in $\Ker (T+1)$:
\begin{eqnarray}
&&u_l\otimes u_m\ +\  u_m\otimes u_l \hspace{2in}
\quad\hbox{if $l\equiv m\bmod n$},
\label{ker-gen-1} \\
&&u_l\otimes u_m\ +\  qu_m\otimes u_l
\ +\  u_{m-i}\otimes u_{l+i}\ +\  q u_{l+i}\otimes u_{m-i}\nonumber\\
&&\phantom{**********}
\hspace{1.2in} \hbox{if $m-l\equiv i\bmod n$ and $0<i<n$}.
\label{ker-gen-2}
\end{eqnarray}

Let $V(z) \wedge_q V(z)$ denote the quotient $(V(z) \otimes V(z))/
\Ker (T+1)$, and let $u_l\wedge_q u_m$ denote
the image of $u_m\otimes u_l$ under the
quotient map. The space
$V(z) \wedge_q V(z)$ will be called {\em $q$-wedge space}
and its elements {\em $q$-wedges}.
{}From now on, we will write $\wedge$ instead of
$\wedge_q$, but it should be understood that {\em all wedges appearing
in this paper are really $q$-wedges.}

The relations
(\ref{ker-gen-1}) and (\ref{ker-gen-2}) can be understood
as normal ordering rules, i.e., as
prescriptions for writing a $q$-wedge whose left entry is smaller than
its right in the ordering (\ref{ordering}) (i.e., $u_l\wedge u_m$
such that $l<m$) as a linear combination of
{\em normally ordered} $q$-wedges whose left entries
are larger than their right.
For $l\equiv m\bmod n$, the rule is simply
\bea
u_l\wedge u_m=-u_m\wedge u_l.
\label{wedge-rel-1}
\ena
In order to give the rule for the case
$m-l\equiv i\bmod n$ and $0<i<n$, let us extract from
(\ref{ordering}) the subsequence
\begin{equation}
\cdots>u_m>u_{m-i}>u_{m-n}>u_{m-n-i}>\cdots
>u_{l+n+i}>u_{l+n}>u_{l+i}>u_l>\cdots
\label{subsequence}
\end{equation}
The rule is
\begin{eqnarray}
u_l\wedge u_m &=& -qu_m\wedge u_l
\ +\ (q^2-1) (u_{m-i}\wedge u_{l+i}
\ -\ qu_{m-n}\wedge u_{l+n} \label{wedge-rel-2} \\
&& \hspace{1.8in} +\ q^2u_{m-n-i}\wedge u_{l+n+i}\ +\cdots ).\nonumber
\end{eqnarray}
Each wedge in the sum is obtained from the
one before it by moving its left component one to the right
in the sequence (\ref{subsequence}),
while simultaneously moving the right component one to
the left.  The sum continues as long as we get normally ordered wedges.

As in the case $N=2$, define $N$-fold
$q$-wedge space $\bigwedge^N V(z)$
to be the quotient $V(z)^{\otimes N}/{\rm Ker}A^{(N)}$.
The next chain of arguments will enable us to conclude that
{\em $\bigwedge^N V(z)$
is equal to the quotient of
$V(z)^{\otimes N}$ by the relations
(\ref{wedge-rel-1}) and (\ref{wedge-rel-2}) in each pair of adjacent
factors}.

The notion of normal ordering was introduced with the motivation that
using the relations (\ref{wedge-rel-1}) and (\ref{wedge-rel-2}),
any element of
$\bigwedge^N V(z)$ can be written as a sum of normally ordered
$q$-wedges (i.e., the terms decrease strictly from left to right with
respect to the ordering (\ref{ordering})).  This means, at the very
least, that the normally ordered $q$-wedges span $\bigwedge^N V(z)$,
but in fact more is true:
\begin{prop}\label{INDP}
The elements
\begin{equation}
u_{m_1}\wedge u_{m_2}\wedge\cdots\wedge u_{m_k},
\label{NOBAS}
\end{equation}
where
$m_1>m_2>\cdots>m_k$, form a basis for $\bigwedge^N V(z)$.
\end{prop}
\pf
In light of the previous discussion,
it remains to show
the normally ordered $q$-wedges are linearly independent.
Because of Proposition \ref{TSPLT}, we have an isomorphism
\[
\bigwedge \! ^N V(z)\simeq\Im A^{(N)}\subset V(z)^{\otimes N}.
\]
Hence the proposition reduces to the linear independence of
the vectors
$$(u_{m_1}\otimes u_{m_2}\otimes\cdots\otimes u_{m_k})A^{(N)}$$
with the $m_i$ strictly decreasing.
This is easily seen by specializing at $q=1$.
\qed

We can conclude from
the independence of the normally ordered $q$-wedges
that the relations (\ref{wedge-rel-1}) and (\ref{wedge-rel-2})
(applied in adjacent factors) generate the entire kernel of $A^{(N)}$.
In other words, they are precisely
the complete set of relations for the $q$-wedge
product that we have been seeking.

%

\subsection{The thermodynamic limit}
Now consider an infinite tensor product $V(z) \otimes V(z) \otimes
V(z) \otimes \cdots$.  (Physicists call this the ``thermodynamic
limit.'')  An infinite iteration of (\ref{coprod-1})-(\ref{coprod-3}) gives
rise to a formal action of $\qslhat$ in this tensor product.  The action
is only formal because when $E_i$ or $F_i$ from $\qslhat$ acts on an element
of the infinite tensor product, the result is typically an infinite sum.
Consequently, it is often not possible to compose two elements of
$\qslhat$. For example, consider
\bea
u_{(m)} = u_m \otimes u_{m-1} \otimes u_{m-2} \otimes \cdots\,.
\label{GST}
\ena
The action of $F_i$ on $u_{(m)}$ produces infinitely many terms:
\begin{equation}
F_iu_{(m)}=q^c\sum
_{\scriptstyle{k\equiv i\,\bmod\,n}\atop\scriptstyle{k\le m}}
u_m\otimes u_{m-1}\otimes\cdots\otimes u_{k+1}
\otimes u_{k+1}\otimes u_{k-1}\otimes\cdots,
\label{FIACT}
\end{equation}
where $c=0$ if $m \equiv i \bmod n$, and $1$ otherwise.
If we apply $E_i$ to the right hand side of (\ref{FIACT}),
all the terms contribute to $u_{(m)}$ because
$E_iu_{k+1}=u_{k}$ for $k\equiv i \bmod n$.
Therefore $E_iF_i$ diverges.

The affine Hecke algebra action behaves better.  The formulas given by
(\ref{hecke-action}) define an action
of the infinite affine Hecke algebra $\hinfhat$ (generated
by $T_i$ and $y_i^{\pm 1}$, $i=1,2,3,\ldots$ with the above relations) on
$V(z) \otimes V(z) \otimes V(z) \otimes \cdots$.  This
action is well-defined because each $T_i$ acts only in a pair
of adjacent factors.
The action of $\hinfhat$
in the thermodynamic limit commutes with the formal action of $\qslhat$.

Let $U_{(m)}$ denote the linear span of all pure tensors that coincide with
$u_{(m)}$ given by (\ref{GST}) after finitely many factors.
In other words, $U_{(m)}$ is spanned by tensors of the form
\[
u_{m_1}\otimes u_{m_2}\otimes u_{m_3}\otimes\cdots
\]
where $m_k=m-k+1$ for $k>\hskip-2pt>1$.

Define $F_{(m)}$ to be the quotient of $U_{(m)}$ by the space
$\sum_i \Ker (T_i+1)$, or, equivalently, by the relations
(\ref{wedge-rel-1}) and (\ref{wedge-rel-2}) in each pair
of adjacent factors.  The spaces
$F_{(m)}$ will be called {\em $q$-deformed Fock spaces}, or
semi-infinite $q$-wedge spaces. Corresponding to (\ref{GST}), we set
\bea
\vac{m} = u_m \wedge u_{m-1} \wedge u_{m-2}\wedge\cdots.
\label{WGST}
\ena

\bigskip
\noindent
{\bf Remark} \hspace{2mm}
In \cite{Eug}, semi-infinite wedges were defined as completely
$q$-antisymmetrized `ideal' elements of $U_{(m)}$,
using the antisymmetrization
operator $\sum_{\sigma \in S_{\infty}} T_{\sigma}$.  ($S_{\infty}$ is
the group of bijections $\mbox{\bf Z}^+ \to \mbox{\bf Z}^+$ fixing
all but finitely many elements; equivalently, it is the group
generated by adjacent transpositions $\sigma_1, \sigma_2, \sigma_3,
\ldots \ $. The elements
$\{T_{\sigma}\}_{\sigma \in S_{\infty}}$ form a basis
for the infinite Hecke algebra $H_{\infty}(q^2)$ inside $\hinfhat$.)
We prefer to define semi-infinite wedge space as a quotient
by certain relations so as not to have to work with infinite sums.

\medskip

%
%
%
Because each $T_i$ is an intertwiner,
the action of $\qslhat$
on $V(z)^{\otimes N}$ given by the coproduct (iterated $N-1$ times)
factors through to the quotient space
$\bigwedge^N V(z)=V(z)^{\otimes N}/\Ker A^{(N)}$.
We will show that the formal action of $\qslhat$
on the infinite tensor space $U_{(m)}$ induces a
genuine action on $F_{(m)}$.

For each vector $v\in F_{(m)}$ we have a decomposition
of the form
\begin{equation}
v=v^{(N)}\wedge \vac{m-N},
\qquad v^{(N)}\in\bigwedge \! ^N V(z)
\label{DCPN}
\end{equation}
for a sufficiently large $N$.
Therefore, if we determine the action of the Chevalley generators
on $\vac{m}$ for all $m\in\Z$, the coproduct
(\ref{coprod-1}-\ref{coprod-3})
gives the action on general vectors.

We define the action of $E_i$, $F_i$, $K_i$ $(i=0,1,\cdots,n-1)$
on $\vac{m}$ to mirror their formal action in $U_{(m)}$,
as follows:
\begin{eqnarray}
E_i\vac{m}&=&0,\\
F_i\vac{m}&=&\cases{
u_{m+1}\wedge u_{m-1}\wedge u_{m-2}\wedge\cdots&if $i\equiv m\bmod n$;\cr
0&otherwise,\cr}\\
K_i\vac{m}&=&\cases{
q\vac{m}&if $i\equiv m\bmod n$;\cr
\vac{m}&otherwise.\cr}
\end{eqnarray}
Noting that $u_l\wedge u_l=0$, we can show the well-definedness of this
action. The action of $d$ is also consistently defined by fixing
the degree of $\vac{0}$ to be zero.

\begin{prop}
The $\qslhat$-module $F_{(m)}$ is isomorphic to the $q$-deformed Fock space
constructed by Hayashi in \cite{Hayashi} and Misra-Miwa in \cite{Misra-Miwa}.
\end{prop}

\pf
{}From Proposition \ref{INDP} we see that the vectors
\[
u_{m_1}\wedge u_{m_2}\wedge u_{m_3}\wedge\cdots
\]
where $m_1>m_2>m_3>\cdots$ and $m_k=m-k+1$
for $k>\hskip-2pt>1$,
constitute a basis of $F_{(m)}$. There is an evident one-to-one
correspondence between these vectors and the colored Young diagrams in
\cite{Misra-Miwa}, and this
correspondence is equivariant with respect to the $\upqn$-actions.
\qed

Set
\[
\goth{h}=\C H_0\oplus\C H_1\oplus\cdots\C H_{n-1}\oplus\C d,
\]
where $K_i=q^{H_i}$.
Let $\Lambda_i\in \goth h^*$ $(0\le i\le n-1)$ be the fundamental weights,
i.e., $\Lambda_i(H_j)=\delta_{i,j}$ and $\Lambda_i(d)=0$.
We define $\Lambda_m$ for $m\in\Z$ by requiring
$\Lambda_m=\Lambda_{m-1}+\hbox{\rm wt\,}u_m$.
(If $m_1 \equiv m_2 \bmod n$, then $\Lambda_{m_1}$ and $\Lambda_{m_2}$
are the same apart from the action of $d$.)
Then $\vac{m}$ has weight $\Lambda_m$
and the weights of $F_{(m)}$ belong to
$\Lambda_m+\sum_{i=0}^{n-1}\Z_{\le 0}\,\alpha_i$.
The highest weight vector $\vac{m}\in F_{(m)}$
generates the irreducible highest weight module $V_{\Lambda_m}$
with highest weight $\Lambda_m$.

It is important to remark that
$\vac{m}$ is not the only highest weight vector in $F_{(m)}$.
For example, if $n=2$, then $(F_0 F_1 - q F_1 F_0) \vac{0}$ is a
highest weight vector that lies in the same Fock space as $\vac{0}$.
The goal of this paper is to describe the
decomposition of $F_{(m)}$ (which is completely reducible) as a
$\uqn$-module.

\section{Heisenberg algebra and decomposition of $F_{(m)}$}
When $q=1$, $F_{(m)}$ reduces to the ordinary infinite wedge space,
and its decomposition as an $\slhat$-module is known.
(See, for example, \cite{q=1}.)
The decomposition comes from a Heisenberg
algebra $H$ acting on wedge space and commuting with the action of
$\slhat$.  In this section, we introduce an
analogous action of the Heisenberg algebra on the $q$-deformed
Fock spaces, and use this action to decompose these spaces as
$\qslhat$-modules.
We also decompose the mapping between Fock spaces that is induced by
the $q$-wedging operator into the product of the vertex operator for
the level-$1$ $\uqn$-modules and that for the Heisenberg algebra.

\subsection{Center of the affine Hecke algebra}
The aim of this section is to define an action of a Heisenberg algebra
$H$ on $F_{(m)}$, which commutes with the action of $\upqn$.
This Heisenberg algebra is the limit $N \rightarrow \infty$ of the
center of the finite affine Hecke algebra $\widehat H_N(q^2)$.
Then the $q$-deformed Fock space $F_{(m)}$,
regarded as a representation of $\upqn\otimes U(H)$, decomposes
into the tensor product
\bea
F_{(m)}\simeq V_{\Lambda_m}\otimes\C[H_-],
\ena
where $\C[H_-]$ is the Fock space of the Heisenberg algebra $H$.

In the previous section we constructed the $q$-deformed Fock space
by starting from the $\uqn$-module $V(z)$. The $\upqn$-action
on $V(z)$ commutes with $y\in\End_{\C}V(z)$ (acting as multiplication
by $z^{-1}$). Using this fact
as a building block, we will construct
elements in the centralizer of the $\upqn$-action on $F_{(m)}$.

For $a\in\Z\backslash\{0\}$ define an operator $B_{a}$ acting
formally in $V(z) \otimes V(z) \otimes V(z) \otimes \cdots$ by
\bea
B_{a}=\sum_{k=1}^\infty y_k^{-a}.
\ena
It is clear
that this formal action commutes with the formal action of $\qslhat$,
and also that it preserves each subspace $U_{(m)}$.
The following can be checked using relations
(\ref{heckegen-6}) and (\ref{heckegen-7}):
\begin{lemma}
The element
\[
B^{(N)}_{a}=\sum^N_{k=1}y_k^{-a}
\]
belongs to the center of $\widehat H_N(q^2)$.
\end{lemma}

In particular, $B_{a}$ commutes with each $T_i$, and therefore
preserves the spaces $\Ker(T_i+1)$.
It therefore preserves their sum $\sum_i\Ker(T_i+1)$,
which means that it acts
on the quotient spaces $F_{(m)}$.  In fact, this is a genuine action
rather than just a formal one: the action of $B_{a}$ on a wedge results
in a finite sum of wedges.  This is because for sufficiently large $k$
we have $y_k^{-a}w=0$, which follows from
\begin{lemma}\label{UDL}
Let $l\le m$. Then, the $q$-wedges
$u_m\wedge u_{m-1}\wedge\cdots\wedge u_{l+1}\wedge u_l\wedge u_m$
and
$u_l\wedge u_m\wedge u_{m-1}\wedge\cdots\wedge u_{l+1}\wedge u_l$
are both equal to zero.
\end{lemma}

\noindent
{\em Proof.} \hspace{2mm}
A straightforward induction using the relations
(\ref{wedge-rel-1}) and (\ref{wedge-rel-2}). \qed

As already mentioned, the action of $B^{(N)}_{a}$ on $V(z)^{\otimes N}$
commutes with the action of $\upqn$. Thus, we conclude
\begin{prop}
The operator $B_{a}$ acts on $F_{(m)}$ and commutes with $\upqn$.
\end{prop}

Next we compute the commutator $[B_{a_1},B_{a_2}]$. For this purpose we need
\begin{lemma}
\label{OPBB}
\bea
[B_{a},z^bv_i]=z^{a+b}v_i.
\ena
\end{lemma}

\noindent
Here we consider $z^bv_i$ as an operator $F_{(m)} \to F_{(m+1)}$
acting as follows: if $v \in F_{(m)}$, then
\[
z^b v_i: v \mapsto z^bv_i\wedge v.
\]
\begin{lemma}
\label{IDL}
Suppose that $\beta\in\End_\C\left(\oplus_{m\in\Z}F_{(m)}\right)$
satisfies the following:
\begin{eqnarray}
&&\hbox{(i)} \qquad \beta F_{(m)}\subset F_{(m)},\label{CN1}\\
&&\hbox{(ii)} \qquad
[d,\beta]=a\beta\quad\hbox{for some $a\in\Z$},\label{CN2}\\
&&\hbox{(iii)} \qquad
[\beta,z^bv_i]=0\quad \hbox{for any $b\in\Z$ and $i\in\{1,2,\ldots,n\}$}
\label{CN3}.
\end{eqnarray}
Then $\beta=\gamma\,\hbox{\rm id}$ for some constant $\gamma$.
\end{lemma}

\pf
We use the decomposition (\ref{DCPN}).
{}From (\ref{CN3}) we have
\[
\beta v=v^{(N)}\wedge\beta \vac{m-N}.
\]
{}From (\ref{CN1}) we have
\[
\beta \vac{m-N}=\sum_k\gamma_ku_{m_{1,k}}\wedge\cdots\wedge u_{m_{M,k}}
\wedge \vac{m-N-M}
\]
for some $M>0$, where $m_{1,k}>\cdots>m_{M,k}>m-N-M$ for all $k$.
Because of (\ref{CN2}), there exists an integer $L$ independent of $N$
such that $m_{1,k}\le m-N+L$.
By Lemma \ref{UDL},
for a sufficiently large $N$, $v^{(N)}\wedge u_{m_{1,k}}=0$
for all $k$ such that $m_{1,k}\ge m-N+1$.
Therefore, we can ignore these terms in $\beta v$.
Suppose that $m_{1,k}\le m-N$ and
$u_{m_{1,k}}\wedge\cdots\wedge u_{m_{M,k}}\wedge \vac{m-N-M}\not=0$.
In this case, Lemma \ref{UDL}
implies $m_{l,k}=m-N-l+1$ $(1\le l\le M)$.
Therefore, $\beta v=\gamma v$, and it is clear that $\gamma$ is independent
of $v$.
\qed

We are now ready to show
\begin{prop}
\begin{eqnarray}
&&\hbox{(i)} \qquad\hbox{If $a_1+a_2\not=0$, then $[B_{a_1},B_{a_2}]=0$.}
\nonumber\\
&&\hbox{(ii)} \qquad
\hbox{$\gamma_a=[B_a,B_{-a}]$ is a non-zero constant.}\nonumber
\end{eqnarray}
\end{prop}
\pf
The basic fact, which we will use to prove (i) and (ii), is that
$[B_{a_1},B_{a_2}]$ is a constant for all $a_1$ and $a_2$.  This will
follow from Lemma \ref{IDL} as soon as we have verified that
the operator $[B_{a_1},B_{a_2}]$ satisfies the three conditions
of the lemma.  It is clear that $[B_{a_1},B_{a_2}]$ preserves $F_{(m)}$.
Since $[d,B_{a}]=aB_{a}$, $[d, [B_{a_1},B_{a_2}]] = (a_1+a_2)
[B_{a_1},B_{a_2}]$, so the second condition is satisfied.
Finally, Lemma \ref{OPBB} implies
\[
\left[[B_{a_1},B_{a_2}],z^bv_i\right]=0,
\]
which is precisely the third condition.

Now to show (i), observe that
if $a_1+a_2\not=0$, then the degree of the operator
$[B_{a_1},B_{a_2}]$ is non-zero. 
Since it is a constant, we must have $[B_{a_1},B_{a_2}]=0$.

To show (ii), note that one can compute
$[B_{a},B_{-a}]\vac{m}$ as a finite sum
by using the formulas (\ref{wedge-rel-1})
and (\ref{wedge-rel-2}).  This implies that
$\gamma_a$ is a polynomial in $q$,
and by specializing it to $q=1$, we can conclude
that $\gamma_a\not=0$. \qed

We will derive an explicit formula for $\gamma_a$
in the next subsection.  For now, let us calculate
$\gamma_1$ as an example.

Using \ref{UDL}, we have
\bea
B_{-1}\vac{0} &=& u_n\wedge \vac{\! - \! 1}
\ + \ u_0\wedge u_{n-1}\wedge \vac{\! - \! 2}
\ + \ u_0\wedge u_{-1}\wedge u_{n-2}\wedge \vac{\! - \! 3} \\
&&\phantom{*******}
+ \cdots+ \ u_0\wedge\cdots\wedge u_{-n+2}\wedge u_1\wedge \vac{\! -\! n}
\nonumber \\
&=& u_n\wedge \vac{\! - \! 1}
\ - \ qu_{n-1}\wedge u_0\wedge \vac{\! -\! 2}
\ + \ (-q)^2u_{n-2}\wedge u_0\wedge u_{-1}\wedge \vac{\! -\! 3} \\
&&\phantom{*******} + \cdots
+ \ (-q)^{n-1}u_1\wedge u_0\wedge\cdots\wedge u_{-n+2}\wedge \vac{\! -\! n}
\nonumber.
\ena
Then, by applying $B_1$ we get
\bea
B_1B_{-1}\vac{0}&=&u_0\wedge \vac{\! -\! 1}
\ - \ qu_{-1}\wedge u_0\wedge \vac{\! -\! 2}
\ + \ (-q)^2u_{-2}\wedge u_0\wedge u_{-1}\wedge \vac{\! -\! 3} \nonumber\\
&&\phantom{*****}
+ \cdots+ \ (-q)^{n-1}u_{-n+1}\wedge u_0\wedge\cdots
\wedge u_{-n+2}\wedge \vac{\! -\! n}
\nonumber\\
&=&(1+q^2+q^4+\cdots+q^{2n-2})\vac{0}.\nonumber
\ena
Therefore, noting that $[B_1, B_{-1}]\vac{0} = B_1 B_{-1}\vac{0}$,
we conclude
\[
\gamma_1={1-q^{2n}\over1-q^2}.
\]

Summing up, in this subsection
we have constructed the action of the Heisenberg
algebra $H$ on $F_{(m)}$.  $H$ is generated by the
operators $B_{a}$ $(a\in\Z\backslash\{0\})$
with the commutation relations
\be
[B_{a},B_{b}]=\delta_{a+b,0}\,\gamma_a, \label{BaBb}
\en
where the constant $\gamma_a$ is yet to be determined.

Let $\C[H_-]$ be the Fock space of $H$, i.e.,
$\C[H_-]=\C[B_{-1},B_{-2},\ldots]$. The element $B_{-a}$
$(a=1,2,\ldots)$ acts on $\C[H_-]$ by multiplication.
The action of $B_{a}$ $(a=1,2,\ldots)$ is given by (\ref{BaBb})
together with the relation
\begin{eqnarray}
&&B_{a} \cdot 1=0\quad\hbox{for $a\ge1$}.
\end{eqnarray}
By computing characters we easily obtain
\begin{prop}
There is an isomorphism
\bea
\iota_m:F_{(m)}\simeq V_{\Lambda_m}\otimes\C[H_-]
\lb{DCM}
\ena
of $\upqn\otimes H$-modules.
\end{prop}
We normalize the isomorphism by requiring
\[
\iota_m(\vac{m})=v_{\Lambda_m}\otimes1.
\]

\subsection{Decomposition of the vertex operator}
We define the vertex operator
\bea
\Omega:V(z)\otimes F_{(m-1)}\rightarrow F_{(m)}
\lb{DECOM}
\ena
by $\Omega(u_a\otimes \omega )=u_a\wedge \omega $.
This is an intertwiner of $\uqn$-modules.

In this section, we decompose the vertex operator $\Omega$
into two parts corresponding to the decomposition \rf{DCM}:
one which acts from $V_{\Lambda_{m-1}}$
to $V_{\Lambda_m}$ and the other which acts on $\C[H_-]$.

The first step in carrying out this decomposition is to transfer
$\Omega$ from the $F_{(m)}$-setting to the $V_{\Lambda_m} \otimes
\C[H_-]$-setting.  To be precise,
define
\bea
\Omega':V(z)\otimes V_{\Lambda_{m-1}}\otimes\C[H_-]
\longrightarrow V_{\Lambda_m}\otimes\C[H_-]
\ena
by requiring that the following diagram commutes:
\bea
\matrix{
V(z)\otimes F_{(m-1)}&
\buildrel{{\rm id}_{V(z)}\otimes\iota_{m-1}}\over\longrightarrow&
V(z)\otimes V_{\Lambda_{m-1}}\otimes\C[H_-]\cr
&&\cr
\Bigg\downarrow\Omega&&\Bigg\downarrow\Omega'\cr
&&\cr \label{DGRM}
F_{(m)}&
\buildrel{\iota_m}\over\longrightarrow&
V_{\Lambda_m}\otimes\C[H_-]\cr}
\ena

We will decompose $\Omega'$ (on the level of generating series)
into one part corresponding to $V_{\Lambda_m}$ and another
part corresponding to $\C [H_-]$.  Given $j\in\Z$, we associate
to $\Omega$ and $\Omega'$ the generating series
\begin{eqnarray*}
\Omega_j(w) &=& \sum_{b\in\Z}\Omega_{j,b}w^{-b}, \\
\Omega'_j(w) &=& \sum_{b\in\Z}\Omega'_{j,b}w^{-b}.
\end{eqnarray*}
Here $\Omega_{j,b}$ is an operator $F_{(m-1)} \longrightarrow F_{(m)}$
whose action is defined by
$$\Omega_{j,b} \cdot \omega = \Omega (u_{j-nb}\otimes \omega).$$
Similarly, $\Omega'_{j,b}$ is an operator $V_{\Lambda_{m-1}} \otimes \C[H_-]
\longrightarrow V_{\Lambda_m} \otimes \C[H_-]$ whose action is given by
$$\Omega'_{j,b}(\omega \otimes f) = \Omega' (u_{j-nb}\otimes \omega
\otimes f).$$

The first element in the decomposition of $\Omega'$ is a
$\upqn$-vertex operator corresponding to $V_{\Lambda_m}$.  It is
given by the following proposition:

\begin{prop}[\cite{vertex}]\label{unicite}
There exists a unique intertwiner of $\uqn$-modules
\bea
\tPhi^*:V(z)\otimes V_{\Lambda_{m-1}}\longrightarrow V_{\Lambda_m}
\ena
such that
$\tPhi^*\Bigl(u_m\otimes v_{\Lambda_{m-1}}\Bigr)=v_{\Lambda_m}$.
\end{prop}

\noindent
The intertwiner $\tPhi^*$ also has a generating series associated
with it.  It is given by
$$\tPhi^*_j(w)=\sum_{b\in\Z}\tPhi^*_{j,b}w^{-b},$$
where $\tPhi^*_{j,b}:V_{\Lambda_{m-1}}\longrightarrow V_{\Lambda_m}$
is given by
\[
\tPhi^*_{j,b}  \cdot \omega =\tPhi^*(u_{j-nb}\otimes \omega ).
\]

\medskip
\noindent
{\bf Remark} \hspace{2mm}
Note that the uniqueness in Proposition
\ref{unicite} is assured by the requirement
that the degree with respect to $d$ is invariant.
A general intertwiner of
$\upqn$-modules
$\Phi:V(z)\otimes V_{\Lambda_{m-1}}\longrightarrow V_{\Lambda_m}$
(not preserving degree) can be written as
\bea
\Phi(z^a v_l\otimes \omega )=
\sum_kc_k
\tPhi^*(z^{a+k} v_l \otimes \omega ).
\ena

The second component in the decomposition
of $\Omega'$ is a vertex operator corresponding
to $H$.  This is given by the generating series
\begin{equation}
\Xi(w)=
\exp\left( \sum_{b\ge1}{B_{-b}w^b\over\gamma_b}\right)
\exp\left( -\sum_{b\ge1}{B_{b}w^{-b}\over\gamma_b}\right). \label{D}
\end{equation}
Here $\gamma_b$ is given by
\bea
\gamma_b=[B_b,B_{-b}].
\lb{E}
\ena
The vertex operator $\Xi(w)$ is characterized by
the commutation relation
\begin{equation}
[B_a,\Xi(w)]=w^a\Xi(w)\,.\label{xi}
\end{equation}
Let $\Xi_b$ denote the coefficient of $w^{-b}$ in
the expansion of (\ref{D}).

We can now state the main result of this section.
\begin{prop}
\bea
\Omega'_j(w)=\tPhi^*_j(w)\otimes\Xi(w).
\lb{C}
\ena \label{main}
\end{prop}
\pf
First we observe that by Proposition \ref{unicite}
the map $\Omega'$ is uniquely determined by
the following conditions:
\begin{eqnarray}
&&\hbox{(i)}\qquad\Omega'(u_m\otimes v_{\Lambda_{m-1}}\otimes1)
=v_{\Lambda_m}\otimes1,\nonumber\\
&&\hbox{(ii)}\qquad\Omega'\hbox{ is $\uqn$-linear},\nonumber\\
&&\hbox{(iii)}\qquad\Omega'\hbox{ is $H$-linear}.\nonumber
\end{eqnarray}
In (iii) the action of $B_{a}$ on
$V(z)\otimes V_{\Lambda_{m-1}}\otimes \C[H_-]$ is given by
\[
B_{a}(z^bv_i\otimes v\otimes f)=
z^{a+b}v_i\otimes v\otimes f+z^av_i\otimes v\otimes B_{a}f.
\]
%
We will show that the right hand side of
\rf{C} satisfies these conditions, and hence is equal to $\Omega'(w)$.
Condition (i) follows from
\bea
\tPhi^*_{m,b}v_{\Lambda_{m-1}}&=&
\cases{0&if $b>0$;\cr v_{\Lambda_m}&if $b=0$,\cr}
\\
\Xi_b1&=&\cases{0&if $b>0$;\cr1&if $b=0$.\cr}
\ena
Condition (ii) is satisfied automatically.
Condition (iii)
follows from (\ref{xi}).
\qed

The final step is to explicitly calculate
$\gamma_a=[B_a,B_{-a}]$
by comparing two point
functions of the vertex operators appearing in Proposition \ref{main}.
The result we are aiming for is

\begin{prop}
\bea
[B_{a},B_{-a}]=a{1-q^{2na}\over1-q^{2a}}.
\ena  \label{comm-formula}
\end{prop}

\pf The idea of the proof is to calculate and compare the two
point functions of each side of (\ref{eqn:C}), then read off a
formula for $\gamma_a = [B_a,B_a]$.  The right hand side can
be done in each factor separately; the answers are given by
the following two lemmas:
\begin{lemma}
\bea
\langle 1,\Xi(w_1)\Xi(w_2)1 \rangle =
{\rm exp}\left(\sum_{a>0}{(w_2/w_1)^a\over\gamma_a}\right)
\label{comxi}
\ena
\end{lemma}
\begin{lemma}
\bea
\langle v_{\Lambda_{m+1}},\tPhi_{m+1}^*(w_1)
\tPhi_m^*(w_2)v_{\Lambda_{m-1}} \rangle =
{(q^{2n+2}w_2/w_1;q^{2n})_\infty\over(q^{2n}w_2/w_1;q^{2n})_\infty}.
\label{comphi}
\ena
Here $(z,p)_\infty=\prod_{k\ge0}(1-zp^k)$.
\end{lemma}
The former is easy, and the latter is obtained in \cite{DO,JM}.  The
two point function of the right hand side of (\ref{eqn:C}) is then just
the product of (\ref{comxi}) and (\ref{comphi}).

Now on to the left hand side of (\ref{eqn:C}).  By (\ref{DGRM}),
the two point function for $\Omega'(w)$ is the same as the one for
$\Omega(w)$.  The latter is given by the following:
\begin{lemma}
\bea
\langle{m+1}|\Omega_{{m+1}}(w_1)\Omega_{m}(w_2)\vac{m-1}
&=&{1-w_2/w_1\over1-q^2w_2/w_1}.
\label{TWOPT}
\ena
\end{lemma}

\pf
We have
\[
\Omega_{m}(w_2)\vac{m-1}
=\sum_{j=0}^\infty u_{m+nj}\wedge \vac{m-1}w_2^{j}.
\]
Applying $\Omega_{{m+1}}(w_1)$ to this sum, and collecting the terms
whose weight is equal to that of $\vac{m+1}$, we get
\[
\sum_{b=0}^\infty u_{m+1-nb}\wedge u_{m+nb}\wedge \vac{m-1}
w_1^{-b}w_2^{b}.
\]
Using the normal ordering rule gives us
\[
\langle{m+1}|(w_2/w_1)^bu_{m+1-nb}\wedge u_{m+nb}\wedge \vac{m-1}
=\cases{1&if $b=0$,\cr q^{2(b-1)}(q^2-1)(w_2/w_1)^b&if $b>0$.\cr}
\]
Summing up for $b$, we obtain (\ref{TWOPT}).
\qed

\medskip

The proof of Proposition \ref{comm-formula} is just a matter of putting
together these three lemmas.  Setting the left and right two point
functions equal to each other and cancelling,
we obtain
\begin{eqnarray}
{\rm exp}\left(\sum_{a>0}{(w_2/w_1)^a\over\gamma_a}\right)
=\prod_{a\ge0}{1-q^{2na}w_2/w_1\over1-q^{2+2na}w_2/w_1}
\end{eqnarray}
A comparison of coefficients results in the asserted formula for
$\gamma_a$.  \qed

\vspace{6mm}
\noindent
{\bf Acknowledgements} \hspace{2mm}
We are grateful to Michio Jimbo and Nicolai Reshetikhin for
helpful discussions.  E.S. also thanks everyone at RIMS for
their excellent hospitality during his visit in May-June of 1995,
when this work was being completed.

\vspace{4mm}
\noindent
{\sc
M.K. and T.M.: Research Institute for Mathematical Sciences, \\
Kyoto University, Sakyo-ku, Kyoto 606, Japan}

\noindent
{\tt masaki@kurims.kyoto-u.ac.jp, miwa@kurims.kyoto-u.ac.jp}

\medskip
\noindent
{\sc
E.S.: Department of Mathematics, UC
Berkeley, Berkeley, CA 94720}

\noindent
{\tt stern@math.berkeley.edu}


\begin{thebibliography}{9}

\bibitem{q=1} E. Date, M. Jimbo, M. Kashiwara, and T. Miwa,
{\em Transformation groups for soliton equations}, in
``Non-linear Integrable Systems --- Classical Theory and
Quantum Theory,'' Proceedings of RIMS Symposium, World
Scientific, 1983, 39-119.

\bibitem{vertex} E. Date, M. Jimbo, and M. Okado,
{\em Crystal bases and $q$ vertex operators},
Commun. Math. Phys. {\bf 155} (1993), 47-69.

\bibitem{DO} E. Date and M. Okado,
{\em Calculation of excitation spectra of the spin model related with
the vector representation of the quantized affine algebra of type
$A^{(1)}_n$.}, Int. J. Mod. Phys. A, 9(1994), 399-417.

\bibitem{Drinfeld} V.G. Drinfeld, {\em Hopf algebra and the
quantum Yang-Baxter equation}, Dokl. Acad. Nauk USSR {\bf 283}
(1985), 1060-1064.

\bibitem{GRV} V. Ginzburg, N. Reshetikhin, and E. Vasserot, {\em
Quantum groups and flag varieties}, Contemp. Math. {\bf 175} (1994),
101-130.

\bibitem{Hayashi} T. Hayashi, {\em $Q$-analogues of Clifford and
Weyl algebras --- spinor and oscillator representations of quantum
enveloping algebras}, Commun. Math. Phys. {\bf 127} (1990), 129-144.

\bibitem{Jimbo} M. Jimbo, {\em A $q$-difference analogue of
$U_q(\goth{g})$ and Yang-Baxter equation}, Lett. Math. Phys. {\bf 10}
(1985), 63-69.

\bibitem{JM} M. Jimbo and T. Miwa, ``Algebraic Analysis of Solvable
Lattice Models,'' CBMS Regional Conference Series in Mathematics.
Vol. 85, A.M.S., Providence, 1995.

\bibitem{Kac-Raina} V. Kac and A. Raina, ``Bombay Lecture on
Highest Weight Representations of Infinite Dimensional Lie
Algebras,'' World Scientific, Singapore, 1987.


\bibitem{Misra-Miwa} K. Misra and T. Miwa, {\em Crystal base for the
basic representation of $U_q (\widehat{\goth{sl}}(n))$}, Commun. Math.
Phys. {\bf 134} (1990), 79-88.

\bibitem{Eug} E. Stern, {\em Semi-infinite wedges and vertex operators},
Internat. Math. Res. Notices {\bf 4} (1995), 201-220.

\end{thebibliography}
\end{document}